\documentclass[11pt]{article}
\usepackage{amssymb,amsmath}
\title{On Finsler surfaces of constant flag curvature with a Killing field
\thanks{%
This work is supported by the National Natural Science Foundation of China 11371032.
Additionally, the first author acknowledges support via DMS-135958
from the United States National Science Foundation
and via a research grant from Duke University.
 }}
\footnotesize{\author{R. L. BRYANT, L. HUANG and X. MO\\
Mathematics Department\\
Duke University, Durham, NC 27708-0230, USA\\
E-mail:  bryant@math.duke.edu\\
School of Mathematical Sciences\\Nankai University, Tianjin
300071, China\\E-mail: huanglb@nankai.edu.cn\\
and\\
Key Laboratory of Pure and Applied Mathematics\\
School of Mathematical Sciences\\ Peking  University,
Beijing 100871, China\\
E-mail:  moxh@pku.edu.cn}
}
\date{}
\begin{document}
\maketitle
\begin{center}{\bf Abstract} \end{center}
\begin{center}
\begin{minipage}{130mm}  
 We study two-dimensional Finsler metrics of constant flag curvature
 and show that such Finsler metrics that admit a Killing field
 can be written in a normal form that depends on two arbitrary functions of one variable.
 Furthermore, we find an approach to calculate these functions for
 spherically symmetric Finsler surfaces of constant flag curvature.
 In particular, we obtain the normal form of the Funk metric on the unit disk $\mathbb{D}^2$.



\end{minipage}
\end{center}

\section{Introduction}

In Riemannian geometry, one has the concept of sectional
curvature. Its analogue in Finsler geometry is called {\em flag
curvature}. A Finsler metric $F$ is said to be of {\em constant (flag) curvature}
if the flag curvature $K=$ constant.
 One of the fundamental problems in Finsler geometry is to
 study Finsler metrics of constant (flag) curvature
because Finsler metrics of constant flag curvature
are the natural generalization of Riemannian metrics of constant sectional curvature.
Recently, great progress has been made in studying Finsler metrics of constant curvature.
The classification of Randers metrics with constant flag curvature has been completed
by D. Bao, C. Robles and Z. Shen [3].
These metrics include the Funk metric on the unit ball and Katok examples [11].
In [19,\,18], X. Mo found many new Finsler metrics of constant flag curvature
by finding Killing fields of generic Bryant metrics and
Mo-Shen-Yang metrics via the navigation problem.

Killing fields on a Finsler manifold $M$ are vector fields
induced by local $1$-parameter group of isometric transformations of $M$.
They are the natural generalization of Killing fields on a Riemannian manifold
and thereby are important in both mathematics and physics.

For instance, Li-Chang-Mo related some Killing fields of Finsler metrics
to the symmetry of very special relativity (VSR for short).
They find that the isometry group of a class of $(\alpha,\,\beta)$-manifold
is the same as the symmetry of VSR [14].
Very special relativity is an interesting approach to investigating
the violation of Lorentz invariance that was developed by Cohen-Glashow [7].
Let $\phi(s)=1+s$. Then $(\alpha,\,\beta)$-manifold
$\left(M,\,\alpha\phi\left(\frac{\beta}{\alpha}\right)\right)$ becomes a Randers manifold.

The main purpose of this paper
is to study Finsler surfaces of constant flag curvature $K$ with a Killing field.
By using moving frame theory, we establish normal forms of such Finsler surfaces
with $K>0$, $K=0$ and $K<0$ respectively (see (4.16), (5.10) and (6.15) below).
In general, the normal form of a class of Finsler metrics
clarifies our understanding of such spaces of Finsler metrics [5].
After noting these normal forms, we obtain the following:

\

\noindent{\bf Theorem 1.1}\,\,
{\em The space of isometry classes of $2$-dimensional Finsler metrics of constant curvature
that admit a Killing field depends on two arbitrary functions of one variable.
}

\

It is worth mentioning that there are many $2$-dimensional Finsler metrics with a non-zero Killing field.
Let $F$ be a Finsler metric on $\mathbb{B}^2(\mu)$, the ball of radius~$\mu$ in~$\mathbb{R}^2$.
$F$ is said to {\em spherically symmetric} if it satisfies
$$
F(Ax,\,Ay)=F(x,\,y)
$$
for all $x\in \mathbb{B}^2(\mu)$, $y\in T_x\mathbb{B}^2(\mu)$ and
$A\in \mathrm{O}(2)$.
Spherically symmetric Finsler metrics $F$ admit the Killing field
$$
X=-x^2\frac{\partial}{\partial x^1}+x^1\frac{\partial}{\partial x^2}\,,
$$
where $x=(x^1,\,x^2)$ (see (7.3) below).
Recently, the study of spherically symmetric
Finsler metrics has attracted considerable attention.
The classification of projective spherically symmetric
Finsler metrics with constant curvature
has just been completed recently by Zhou, Mo-Zhu and Li [22,\, 20,\,13].
The following expression for spherically symmetric
Finsler metrics had been obtained by Huang-Mo [10,\,11] :
$$
F=|y|\phi\left(\frac{|x|^2}2,\, \frac{\langle x,\,y\rangle}{|y|}\right).
$$
This expression motivates us to find an approach to calculate
two functions of one variable in the normal form of spherically symmetric
Finsler surfaces of constant flag curvature (see Section 7).
In particular, we will obtain the normal form
of the Funk metric on the unit disk $\mathbb{D}^2$.

We will determine the space of Finsler surfaces of constant flag curvature
that admit two linearly independent Killing fields in a forthcoming paper.

\section{Preliminaries}


\subsection{The structure equations of a Finsler surface}

Let $(M,\,F)$ be an oriented Finsler surface.
The function $F$ determines and is determined by the set
$$
\Sigma=\left\{(x,\,y)\in TM\mid F(x,\,y)=1\,\right\}
$$
which is known as the unit tangent bundle of $F$.
For each $x\in M$, the intersection $\Sigma_x=\Sigma\cap T_xM$
is the {\em indicatrix}.
Define
$$
\omega_1:=F_{y^j}dx^j.
\eqno(2.1)
$$
Then $\omega_1$ is a differential form on $\Sigma$.
The form $\omega_1$ is known in the calculus of variations as the {\em Hilbert form}.
On $\Sigma$, there exists a canonical coframing $\omega = (\omega_1,\,\omega_2,\,\omega_3)$
satisfying the {\em structure equations}
$$
d\omega_1=-\omega_2\wedge\omega_3,
\eqno(2.2)
$$
$$
d\omega_2=-\omega_3\wedge\omega_1+I\,\omega_3\wedge\omega_2,
\eqno(2.3)
$$
$$
d\omega_3=-K\,\omega_1\wedge\omega_2-J\,\omega_2\wedge\omega_3
\eqno(2.4)
$$
where the functions $I$, $J$ and $K$ are known as the {\em main scalar},
the {\em Landsberg curvature} and the {\em flag curvature} respectively [4,\,16].


Conversely, if $\Sigma$ is a $3$-manifold
endowed with a coframing $\omega = (\omega_1,\omega_2,\omega_3)$
that satisfies the structure equations (2.2--4)
for some functions $I$, $J$, and $K$ on~$\Sigma$
and has the property (which always holds locally)
that there exists a smooth submersion~$\pi:\Sigma\to M$,
where~$M$ is a surface, whose fibers are the integral curves of~$\omega_1$ and~$\omega_2$,
then there is a unique immersion~$\iota:\Sigma\to TM$ compatible with~$\pi$
that realizes~$\Sigma$  as the unit sphere bundle
of a locally defined Finsler structure~$F$ on~$M$
in such a way that the given coframing is the canonical coframing induced on~$\Sigma$
by the (local) Finsler structure~$F$.  In this way, one has a local equivalence
between Finsler surfaces and $3$-manifolds~$\Sigma$
endowed with a coframing~$\omega$ that satisfies (2.2--4).

\subsection{The Bianchi identities}


Differentiating (2.3) and using (2.4), (2.2) and (2.3),
one deduces
$$
(J\,\omega_1-dI)\wedge\omega_2\wedge\omega_3=0.
\eqno(2.5)
$$
Put
$$
dI:=I_1\,\omega_1+I_2\,\omega_2+I_3\,\omega_3.
\eqno(2.6)
$$
From (2.5) and (2.6), we have
$0=(J-I_1)\,\omega_1\wedge\omega_2\wedge\omega_3.$
It follows that $I_1=J.$
Thus we obtain the following Bianchi identity
$$
dI:=J\,\omega_1+I_2\,\omega_2+I_3\,\omega_3.
\eqno(2.7)
$$
Differentiating (2.4) and using the structure equations (2.4), (2.2) and (2.3),
one obtains
$$
-(IK+J_1+K_3)\,\omega_1\wedge\omega_2\wedge\omega_3=0
\eqno(2.8)
$$
where
$$
dJ:=J_1\,\omega_1+J_2\,\omega_2+J_3\,\omega_3,\quad dK:=K_1\,\omega_1+K_2\,\omega_2+K_3\,\omega_3.
\eqno(2.9)
$$
It follows that $J_1=-IK-K_3$. Assume that the Finsler surface $(M,\,F)$ has constant flag curvature.
Thus
$J_1=-IK.$
Substituting this into (2.9) yields
$$
dJ:=-KI\,\omega_1+J_2\,\omega_2+J_3\,\omega_3.
\eqno(2.10)
$$


\section{Finsler metrics with a Killing field}

Assume
that the Finsler surface $(M,\,F)$ has constant flag curvature~$K$ and
that $(M,\,F)$ admits a non-zero Killing field~$X$.
Then its flow $\varphi_t$ is an isometry on~$(M,\,F)$,
i.e. $\check{\varphi_t}F=F$ where $\check{\varphi_t}$
is the flow on $TM$ defined by
$\check{\varphi_t}(x,\,y):=\left(\varphi_t(x),\,\varphi_{t*}(y)\right).$
Note that $\omega_j$ are globally defined on $\Sigma$ and $\varphi_t$ preserves the orientation and $\Sigma$.
It is easy to see that $\check{\varphi_t}^*\omega_j=\omega_j,\quad j=1,\,2,\,3.$
Thus, we have
$$
\mathcal{L}_{\hat{X}}\omega_j
  =\lim_{t\to 0}\frac 1t \left(\omega_j-\check{\varphi_t}^*\omega_j\right)=0,\quad j=1,\,2,\,3
\eqno(3.1)
$$
where $\hat{X}$ is the natural lift of $X$ to $\Sigma$.
Equation (3.1) tells us that $\hat{X}$ is a symmetry vector field,
that is, a nonzero field whose flow preserves $\omega_j$. It follows that
$$
d(\iota_{\hat{X}}\omega_j)+\iota_{\hat{X}}(d\omega_j)=(d\circ\iota_{\hat{X}}+\iota_{\hat{X}}\circ d)\omega_j=\mathcal{L}_{\hat{X}}\omega_j=0
\eqno(3.2)
$$
where $\iota_{\hat{X}}$ is the interior product with respect to $\hat{X}$.

Write 
$$
a_j=\omega_j(\hat{X})=\iota_{\hat{X}}\omega_j.
\eqno(3.3)
$$
Using (3.2) and (3.3), we have
$$
da_j=d(\iota_{\hat{X}}\omega_j)=-\iota_{\hat{X}}d\omega_j,\quad j=1,\,2,\,3.
\eqno(3.4)
$$
Applying the structure equations, we have
$$
da_1
=
-\iota_{\hat{X}}\left(-\omega_2\wedge\omega_3\right)\\
=
\left(\iota_{\hat{X}}\omega_2\right)\omega_3-\left(\iota_{\hat{X}}\omega_3\right)\omega_2=a_2\omega_3-a_3\omega_2,
\eqno(3.5)
$$
$$
\begin{array}{ccl}
da_2
&=&
-\iota_{\hat{X}}\left(-\omega_3\wedge\omega_1-I\omega_2\wedge\omega_3\right)\\
&=&
\left(\iota_{\hat{X}}\omega_3\right)\omega_1-\left(\iota_{\hat{X}}\omega_1\right)\omega_3+I\iota_{\hat{X}}\left(\omega_2\wedge\omega_3\right)\\
&=&
a_3\omega_1-a_1\omega_3+Ida_1,
\end{array}
\eqno(3.6)
$$
$$
\begin{array}{ccl}
da_3
&=&
-\iota_{\hat{X}}\left(-K\omega_1\wedge\omega_2-J\omega_2\wedge\omega_3\right)\\
&=&
K\left[\left(\iota_{\hat{X}}\omega_1\right)\omega_2-\left(\iota_{\hat{X}}\omega_2\right)\omega_1\right]+J\iota_{\hat{X}}\left(\omega_2\wedge\omega_3\right)\\
&=&
K(a_1\omega_2-a_2\omega_1)+Jda_1.
\end{array}
\eqno(3.7)
$$
By using (3.7), (3.5) and (3.6), we get
$$
\begin{array}{ccl}
\frac 12 d\left(K a_2^2+a_3^2\right)
&=&
Ka_2da_2+a_3da_3\\
&=&
K a_2\left(a_3\omega_1-a_1\omega_3+Ida_1\right)\\
&&+a_3\left[K(a_1\omega_2-a_2\omega_1)+Jda_1\right]\\
&=&
\left(KIa_2+Ja_3\right)da_1-Ka_1\left(a_2\omega_3-a_3\omega_2\right)\\
&=&
\left(KIa_2+Ja_3-Ka_1\right)da_1.
\end{array}
\eqno(3.8)
$$
It follows that $K a_2^2+a_3^2$ is constant on the level sets of~$a_1$.

Assuming that $da_1$ is nonvanishing,
and that the level sets of~$a_1$ on~$\Sigma$ are connected,
one can write
$$
K a_2^2+a_3^2 = 2f(a_1)
\eqno(3.9)
$$
where $f:a_1(\Sigma)\to\mathbb{R}$ is a smooth function
Differentiating (3.9) and using (3.8), we have
$$
2\,f'(a_1)\,da_1  
=d\bigl(2f(a_1)\bigr)
=d\left(K a_2^2+a_3^2\right)=
2\left(KIa_2+Ja_3-Ka_1\right)da_1.
$$
It follows that
$$
KIa_2+Ja_3=f'(a_1)+a_1K
\eqno(3.10)
$$
on $\left\{p\in \Sigma\mid da_1(p)\neq 0\right\}$.
(One can also treat the case in which $a_1$ is constant.)


Because the flow of $\hat X$ preserves the $\omega_i$,
it must also preserve $I$ and $J$.
Thus,


$$
\mathcal{L}_{\hat{X}}I=0
\eqno(3.11)
$$
and
$$
\mathcal{L}_{\hat{X}}J=0.
\eqno(3.12)
$$
By (2.6) and (3.11), we get
$$
\begin{array}{ccl}
a_1J+a_2I_2+a_3I_3
&=&
\left(\iota_{\hat{X}}\omega_1\right)J+\left(\iota_{\hat{X}}\omega_2\right)I_2
+\left(\iota_{\hat{X}}\omega_3\right)I_3\\
&=&
\left(\iota_{\hat{X}}J\right)\omega_1+J\iota_{\hat{X}}\omega_1
+\left(\iota_{\hat{X}}I_2\right)\omega_2+I_2\iota_{\hat{X}}\omega_2\\
&&
+\left(\iota_{\hat{X}}I_3\right)\omega_3+I_3\iota_{\hat{X}}\omega_3\\
&=&
\iota_{\hat{X}}\left(J\omega_1+I_2\omega_2+I_3\omega_3\right)\\
&=&
\iota_{\hat{X}} dI
=(d\circ\iota_{\hat{X}}+\iota_{\hat{X}}\circ d)I=\mathcal{L}_{\hat{X}}I=0.
\end{array}
\eqno(3.13)
$$
Similarly, (2.9) and (3.12) imply that
$-a_1KI+a_2J_2+a_3J_3=0.$
Together with (3.13), (3.7), (3.5) and (3.6), we obtain
$$
\begin{array}{ccl}
d(a_2J-a_3I)
&=&
J\,da_2+a_2\,dJ-I\,da_3-a_3\,dI\\
&=&
J\left(a_3\omega_1-a_1\omega_3+Ida_1\right)+a_2\left(-KI\omega_1+J_2\omega_2+J_3\omega_3\right)\\
&&
-I\left(Ka_1\omega_2-Ka_2\omega_1+Jda_1\right)-a_3\left(J\omega_1+I_2\omega_2+I_3\omega_3\right)\\
&=&
(a_2J_2-Ka_1I-a_3I_2)\omega_2+(a_2J_3-a_1J-a_3I_3)\omega_3\\
&=&
-a_3(J_3+I_2)\omega_2+a_2(J_3+I_2)\omega_3=(J_3+I_2)da_1.
\end{array}
$$
It follows that
$$
a_2\,J-a_3\,I=g(a_1)
\eqno(3.14)
$$
for some function $g$,
again, assuming that $da_1$ is nonvanishing and that the level sets of $a_1$ are connected.

It turns out to be convenient to split
the further discussion into cases according to whether $K>0,\, K=0$ and $K<0$.
Moreover, by scaling, one can reduce to the cases $K=1,\, K=0$ and $K=-1$.

To simplify notation, in the following sections, we shall abbreviate $a_1$ as $a$.
We will also assume that $da$ is nonvanishing and that the level sets of $a$ are connected.

\section{$K=1$}

In this section, we are going to study Finsler surfaces with constant flag curvature $K=1$.
In this case, $a_2^2+a_3^2$ is a function of $a$ $(:=a_1)$ by (3.9).
Without loss of generality, we can assume that $(a_2,\,a_3)\neq (0,\,0)$.
Let
$$
a_2^2+a_3^2=u(a)^2
\eqno(4.1)
$$
where $u$ is a positive function on $a(\Sigma)\subset\mathbb{R}$.

Write
$$
a_2=u(a)\,\sin t,\quad a_3=u(a)\,\cos t
\eqno(4.2)
$$
where $t:\Sigma\to \mathbb{R}$. It follows from (3.10) that
$$
a_2I+a_3J=u(a)\,u'(a) + a.
\eqno(4.3)
$$
We rewrite (3.14) as follows
$$
a_2J-a_3I=u(a)^2v(a)
\eqno(4.4)
$$
where $v(a):=\frac{g(a)}{u(a)^2}$.
For notational simplicity in what follows,
we will write $u$, $u'$, or $v$ instead of $u(a)$, $u'(a)$, or $v(a)$.

Solving (4.4) and (4.3) and then using (4.2), we obtain
$$
I=\left[u'+\frac{a}{u}\right]\,\sin t-uv\,\cos t,
\eqno(4.5)
$$
$$
J=\left[u'+\frac{a}{u}\right]\,\cos t+uv\,\sin t.
\eqno(4.6)
$$
By (4.2), we have
$$
da_2=u'\,\sin t \,da+u\,\cos t\,dt
\eqno(4.7)
$$
and
$$
da_3=u'\,\cos t\,da-u\,\sin t\,dt.
\eqno(4.8)
$$
Plugging (4.2) into (3.5) yields
$$
da=u\left(\sin t\,\omega_3-\cos t\,\omega_2\right).
\eqno(4.9)
$$
By substituting (4.5) into (3.6) and using (4.7), we obtain
$$
\begin{array}{ccl}
u'\,\sin t\,da+u\,\cos t\,dt
&=&
da_2\\
&=&
a_3\,\omega_1-a\,\omega_3+I\,da\\
&=&
a_3\,\omega_1-a\,\omega_3+\left[\left(u'+{\displaystyle\frac a{u}}\right)\,\sin t-uv\,\cos t \right]\,da.
\end{array}
$$
Together with (4.9), we have
$$
u\,\cos t\,dt
=u\,\cos t\left[\omega_1-\left(\frac a{u} \,\sin t-uv\,\cos t\right)\omega_2
    -\left(\frac a{u} \,\cos t+uv\,\sin t\right)\omega_3\right].
$$
Except where $\cos t=0$, we get
$$
dt=\omega_1-\left(\frac au\,\sin t-uv\,\cos t\right)\omega_2
      -\left(\frac au\,\cos t+uv\,\sin t\right)\omega_3.
\eqno(4.10)
$$
Similarly, we obtain that (4.10) holds when $\sin t\neq 0$
using (4.6), (3.6), (4.8) and (4.9). Hence (4.10) holds on $\Sigma$.

Let
$$
\alpha=\frac{\theta}{u}
\eqno(4.11)
$$
where
$$
\theta=\theta_1+\theta_2
\eqno(4.12)
$$
where
$$
\theta_1=(\sin t)\,\omega_2,\quad \theta_2=(\cos t)\,\omega_3.\quad
\eqno(4.13)
$$
Using (2.2), (4.5) and (4.10) we get
$$
d\theta_1= (\cos t)\,\omega_1\wedge\omega_2-(\sin t)\,\omega_3\wedge\omega_1
+\left(\frac au\,\cos 2t-u'\sin^2t+uv\sin 2t\right)\omega_2\wedge\omega_3.
$$
By (2.3), (4.5) and (4.10), we see that
$$
d\theta_2=-(\cos t)\,\omega_1\wedge\omega_2+(\sin t)\,\omega_3\wedge\omega_1
-\left(\frac au\,\cos 2t+u'\cos^2t+uv\sin 2t\right)\omega_2\wedge\omega_3.
$$
Thus, we have
$$
d\theta=-u'\,\omega_2\wedge\omega_3
\eqno(4.14)
$$
from which, together with (4.9), we obtain that the $1$-form $\alpha$ is closed.
It follows that $\alpha$ is locally an exact differential form,
i.e., locally there exists a function $b$ such that $\alpha=db$.
Using (4.11), (4.12) and (4.13), we get
$$
db=\frac 1{u(a)}\left(\sin t\,\omega_2+\cos t\,\omega_3\right).\eqno(4.15)
$$
Taking this together with (4.9) and (4.10), we obtain
$$
\begin{array}{ccl}
da\wedge db\wedge dt
&=&
u\left(\sin t\,\omega_3+\cos t\,\omega_2\right)\wedge
     \frac 1{u}\left(\sin t\,\omega_2+\cos t\,\omega_3\right)\wedge dt\\
&=&
\omega_3\wedge\omega_2\wedge dt\\
&=&
\omega_3\wedge\omega_2\wedge\omega_1=-dV_{\Sigma}\neq 0.
\end{array}
$$
Hence $(a,\,b,\,t)$ is a local coordinate system on $\Sigma$.
Furthermore, we can solve for the $\omega_i$ in the form
$$
\left(
\begin{array}{c}
  \omega_1 \\
  \omega_2 \\
  \omega_3
\end{array}
\right)
=
\left(
\begin{array}{ccc}
  1& v&a \\
  0&           -(\cos t)/u & u\,\sin t\\
  0& \phantom{-}(\sin t)/u & u\,\cos t
  \end{array}
\right)
\left(
\begin{array}{c}
  dt\\
  da\\
  db
\end{array}
\right)
\eqno(4.16)
$$
where we have used (4.9), (4.10) and (4.15).
We say that
(4.16) is the {\em normal form} for a Finsler surface of constant flag curvature $1$
that admits a Killing field.  

Conversely, regarding $u = u(a)>0$ and $v = v(a)$ as arbitrary functions of $a$,
the above coframing $\omega_i$ satisfies the structure equations of an oriented
Finsler surface with $K\equiv1$ and admitting a symmetry vector field.  Note that
it depends on two arbitrary functions of one variable.
Thus, we have shown the following:

\

\noindent{\bf Theorem 4.1}\, {\em The space of isometry classes of Finsler metrics
of constant flag curvature $1$ that admit a Killing field
depends on two arbitrary functions of one variable}.

\

Now we investigate the geometric meanings of $b$, $t$ and $a$.
Using (4.2), (4.11), (4.12) and (4.13), one can verify that
$$
db\left(\hat{X}\right)=\frac 1u(a_2\sin t+a_3\cos t)=1,
$$
$$
dt\left(\hat{X}\right)=a-\left(\frac au\sin t-uv\cos t\right)a_2-\left(\frac au\cos t-uv\sin t\right)a_3=0,
$$
$$
da\left(\hat{X}\right)=u\left(a_3\sin t-a_2\cos t\right)=0.
$$
By the above formula, we obtain $\hat{X}=\frac{\partial}{\partial b}$.
It follows that the $b$-curves are the integral curves of $\hat{X}$.

Let $E$ be the Reeb vector field of $F$. Then [8, Proposition 3.2]
$$
\omega_1(E)=1,\quad \omega_2(E)=\omega_3(E)=0.
$$
Together with (4.9), (4.10) and (4.15) we get
$$
da(E)=db(E)=0,\quad dt(E)=1.
$$
It follows that $E=\frac{\partial}{\partial t}$, equivalently,
the $t$-curves are the flows of $E$.
Recall that a curve is a (unit Finslerian speed)
geodesic if its canonical lift in $\Sigma$ is an integral curve  of $E$ [4].
Thus the $t$-curves are the canonical lift of unit geodesics on $(M,\,F)$.

Finally we are going to discuss the geometric meaning of $a$.
In natural coordinates, we have
$$
\hat{X}=v^i\frac{\partial}{\partial x^i}
+y^j\frac{\partial v^i}{\partial x^j}\frac{\partial}{\partial y^i}
\eqno(4.17)
$$
where $X=v^i\frac{\partial}{\partial x^i}$; see [17, 18]. Together with (2.1) we have
$$
a=
a_1
=
\omega_1(\hat{X})
=
\left(F_{y^k}dx^k\right)
\left(v^i\frac{\partial}{\partial x^i}
       +y^j\frac{\partial v^i}{\partial x^j}\frac{\partial}{\partial y^i}\right)
=
\omega_1(X)=\iota_X\omega_1.
$$
It follows that $a$ is the interior product with respect to the Killing field $X$
of the Hilbert form~$\omega_1$ [12, Page 35].

\section{$K=0$}

In this section, we are going to investigate Finsler surfaces with a flag curvature $K=0$.
In this case, $a_3^2$, $a_3J$ and $a_2\,J-a_3\,I$
are functions of $a$ by using (3.9), (3.10) and (3.14).
Let
$$
a_3=u(a)
\eqno(5.1)
$$
where $u(a)$ is a non-zero function on $\Omega:=\{p\in \Sigma \mid a_3(p)\neq 0\}$.
By using (3.9), (3.10) and (5.1), we have
$a_3J=\left[\frac{u(a)^2}2\right]'=u(a)u'(a).$
Together with (5.1) we get
$$
J=\frac{u(a)u'(a)}{a_3}=u'(a)
\eqno(5.2)
$$
on $\Omega$. Write
$$
a_2=u(a)t,\quad a_2J-a_3I=u(a)^2v(a)
\eqno(5.3)
$$
where $t:\Sigma\to \mathbb{R}$. It follows from (5.2) that
$$
I
=
\frac{a_2J-u(a)^2v(a)}{a_3}=
\frac{u(a)tu'(a)-u(a)^2v(a)}{a_3}=u'(a)t-u(a)v(a).
\eqno(5.4)
$$

Again, for simplicity of notation, let us write $u$, $u'$ or $v$
instead of $u(a)$, $u'(a)$, or $v(a)$.

By (5.2) and the first equation of (5.3), we have
$$
da_2=u't\,da+u\,dt
\eqno(5.5)
$$
and
$$
da_3=u'\,da.
\eqno(5.6)
$$
Substituting (5.1) and the first equation of (5.3) into (3.4), we have
$$
da=u\,(t\,\omega_3-\omega_2).
\eqno(5.7)
$$
By substituting (5.4) into (3.6) and using (5.5), we obtain
$$
\begin{array}{ccl}
u't\,da+u\,dt
&=&
da_2\\
&=&
a_3\,\omega_1-a\,\omega_3+I\,da\\
&=&  
a_3\,\omega_1-a\,\omega_3+\left(u't-uv\right)\,da.
\end{array}
$$
Together with (5.7), we have
$u\,dt=u\left[\omega_1+uv\,\omega_2-\left(\frac au+uvt\right)\,\omega_3\right].$
It follows that
$$
dt=\omega_1+uv\,\omega_2-\left(\frac au+uvt\right)\,\omega_3
\eqno(5.8)
$$
on $\Omega$. By using (2.3), (3.5), (5.1) and (5.2), we have
$$
\begin{array}{ccl}
d\left(\frac 1{u}\omega_3\right)
&=&
-\frac{u'}{u^2}\,da\wedge\omega_3+\frac 1u\,d\omega_3\\
&=&
-\frac{u'}{u^2}(a_2\omega_3-a_3\omega_2)\wedge\omega_3-\frac 1u J\,\omega_2\wedge\omega_3\\
&=&
-\frac{u'}{u^2}(-u\omega_2\wedge\omega_3)-\frac {u'}u\omega_2\wedge\omega_3=0.
\end{array}
$$
We get that the $1$-form $\frac 1{u}\,\omega_3$ is closed.
Hence locally there exists a function $b$ such that
$$
db=\frac 1{u}\,\omega_3.
\eqno(5.9)
$$
Taking this together with (5.7) and (5.8) we obtain
$$
\begin{array}{ccl}
da\wedge db\wedge dt
&=&
u\left(t\,\omega_3-\omega_2\right)\wedge \frac 1{u}\,\omega_3\wedge dt\\
&=&
-\omega_2\wedge\omega_3\wedge\omega_1=-dV_{\Sigma}\neq 0.
\end{array}
$$
It follows that $(a,\,b,\,t)$ is a local coordinate system on $\Sigma$.
Furthermore, just as in the case $K\equiv1$, we obtain a normal form for $(M,\,F)$
with flag curvature $K\equiv0$ that admits a Killing field
$$
\left(
\begin{array}{c}
  \omega_1 \\
  \omega_2 \\
  \omega_3
\end{array}
\right)
=
\left(
\begin{array}{ccc}
  1& v&a \\
  0& -1/u & t\,u(a)\\
 0& 0 & u
  \end{array}
\right)
\left(
\begin{array}{c}
  dt\\
  da\\
  db
\end{array}
\right),
\eqno(5.10)
$$
where $u=u(a)>0$ and $v=v(a)$ are arbitrary functions of $a$.
Thus, we have the following result:

\

\noindent{\bf Theorem 5.1}\,
{\em The space of isometry classes of Finsler metrics of constant flag curvature $0$
that admit a Killing field depends on two arbitrary functions of one variable}.

\

Just as in the previous case of $K\equiv1$,
the geometric meanings of $b$, $t$ and $a$ are as follows:
the $b$-curves are the integral curves of $\hat{X}$;
the $t$-curves are the canonical lift of unit geodesics on $(M,\,F)$
and $a$ is the interior product with respect to the Killing field $X$
of the Hilbert form $\omega_1$.

\section{$K=-1$}

Now let us consider Finsler surfaces of constant flag curvature $K=-1$.
In this case, $-a_2^2+a_3^2$, $-a_2I+a_3J$, $a_2J-a_3I$
are functions of $a$ by virtue of (3.9), (3.10) and (3.14).
By the sign and continuity of $-a_2^2+a_3^2$,
we should investigate the following three subcases:
$$
{\rm(i)} \,\, -a_2^2+a_3^2>0,\quad
{\rm(ii)} \,\, -a_2^2+a_3^2\equiv 0, \quad
{\rm(iii)} \,\, -a_2^2+a_3^2<0
$$
For brevity, we only discuss the subcase (i), as the others are similar.
Moreover, we will continue to assume that $a = a_1:\Sigma\to\mathbb{R}$
is a submersion with connected fibers.

Let
$$
-a_2^2+a_3^2=u(a)^2
\eqno(6.1)
$$
where $u$ is a positive function on $a(\Sigma)\subset\mathbb{R}$. Write
$$
a_2=u(a)\,\sinh t,\quad a_3=u(a)\,\cosh t
\eqno(6.2)
$$
where $t:\Sigma\to \mathbb{R}$. By using (3.9) and (3.10) we have
$$
-a_2I+a_3J=\left[\frac 12u(a)^2\right]'+(-1)a=u(a)u'(a)-a.
\eqno(6.3)
$$
By (3.14), we have
$$
a_2J-a_3I=u(a)^2v(a)
\eqno(6.4)
$$
where $v(a):=\frac{g(a)}{u(a)^2}$.
Again, for simplicity, we will write $u$, $u'$ or $v$
for $u(a)$, $u'(a)$, or $v(a)$, respectively.
Solving (6.3) and (6.4) and using (6.2), we obtain
$$
I=\left[u'-\frac{a}{u}\right]\,\sinh t-uv\,\cosh t,
\eqno(6.5)
$$
$$
J=\left[u'-\frac{a}{u}\right]\,\cosh t-uv\,\sinh t.
\eqno(6.6)
$$
By (4.2), we have
$$
da_2=u'\,\sinh t \,da+u\,\cosh t\,dt
\eqno(6.7)
$$
and
$$
da_3=u'\,\cosh t da+u\,\sinh t\,dt.
\eqno(6.8)
$$
Substituting (6.2) into (3.5) yields
$$
da=u\left(\sinh t\,\omega_3-\cosh t\,\omega_2\right).
\eqno(6.9)
$$
By substituting (6.5) into (3.6) and using (6.7) we obtain
$$
\begin{array}{ccl}
u'\,\sinh t \,da+u\,\cosh t\,dt
&=&
da_2\\
&=&
a_3\,\omega_1-a\,\omega_3+I\,da\\
&=&   
a_3\,\omega_1-a\omega_3+\left[(u'-\frac{a}{u})\sinh t-uv\,\cosh t \right]da.
\end{array}
$$
Together with (6.9), we have
$$
u\,\cosh t\,dt=u\,\cosh t\,\left[\omega_1+\left(\frac au \,\sinh t+uv\,\cosh t\right)\,\omega_2
      -\left(\frac au \,\cosh t+uv\,\sinh t\right)\,\omega_3\right]
$$
where we have used $1+\sinh^2t=\cosh^2t$.
Note that $u\,\cosh t>0$. Hence
$$
dt=\omega_1+\left(\frac au \,\sinh t+uv\,\cosh t\right)\,\omega_2
       -\left(\frac au \,\cosh t+uv\,\sinh t\,\right)\,\omega_3.
\eqno(6.10)
$$
Let
$$
\alpha=\frac{\theta}{u}
\eqno(6.11)
$$
where
$$
\theta=\theta_2-\theta_1
\eqno(6.12)
$$
where
$$
\theta_1=(\sinh t)\,\omega_2,\quad \theta_2=(\cosh t)\,\omega_3.\quad
\eqno(6.13)
$$
Using (2.3), (6.5) and (6.10) we get
$$
d\theta_1=(\cosh t)\,\omega_1\wedge\omega_2-(\sinh t)\,\omega_3\wedge\omega_1
+\left(\frac au\,\cosh 2t-u'\sinh^2t+uv\,\sinh 2t\right)\,\omega_2\wedge\omega_3.
$$
By (2.4), (4.6) and (6.10), we see that
$$
d\theta_2=(\cosh t)\,\omega_1\wedge\omega_2-(\sin t)\,\omega_3\wedge\omega_1
+\left(\frac au\,\cosh 2t-u'\cosh^2t+uv\,\sinh 2t\right)\,\omega_2\wedge\omega_3\,.
$$
Thus, we obtain (4.14), where $\theta$ is defined in (6.12).
By using (4.14) and (6.9) we obtain that the $1$-form $\alpha$ is closed.
It follows that locally there exists a function $b$ such that
$\alpha=db$. Using (6.11), (6.12) and (6.13), we get
$$
db=\frac 1{u}\left(\cosh t\,\omega_3-\sinh t\,\omega_2\right).\eqno(6.14)
$$
Together with (6.9) and (6.10) yields
$$
\begin{array}{ccl}
da\wedge db\wedge dt
&=&
u\left(\sinh t\,\omega_3-\cosh t\,\omega_2\right)\wedge
     \frac 1{u}\left(\cosh t\,\omega_3-\sinh t\,\omega_2\right)\wedge dt\\
&=&
(\cosh^2t-\sinh^2t)\omega_3\wedge\omega_2\wedge \omega_1\\
&=&
\omega_3\wedge\omega_2\wedge\omega_1=-dV_{\Sigma}\neq 0.
\end{array}
$$
Hence $(a,\,b,\,t)$ is a local coordinate system on $\Sigma$.
Moreover we can solve for the $\omega_i$ in the form
$$
\left(
\begin{array}{c}
  \omega_1 \\
  \omega_2 \\
  \omega_3
\end{array}
\right)
=
\left(
\begin{array}{ccc}
  1& v&a \\
  0& -(\cosh t)/u & u\,\sinh t\\
  0& -(\sinh t)/u & u\,\cosh t
  \end{array}
\right)
\left(
\begin{array}{c}
  dt\\
  da\\
  db
\end{array}
\right)
\eqno(6.15)
$$
where we have used (6.9), (6.10) and (6.14).
Then (6.15) is the {\em normal form} for a Finsler surface of constant flag curvature
$K\equiv-1$ that admits a Killing field.
It depends on two arbitrary functions of one variable, namely $u=u(a)>0$
and $v=v(a)$. We then have the following:

\

\noindent{\bf Theorem 6.1}\,
{\em The space of isometry classes of Finsler metrics
of constant flag curvature $-1$ which admits a Killing field
depends on two arbitrary functions of one variable}.

\

Again, as in the earlier cases
the geometric meanings of $b$, $t$ and $a$ are as follows:
the $b$-curves are the integral curves of $\hat{X}$;
the $t$-curves are the canonical lift of unit geodesics on $(M,\,F)$
and $a$ is the interior product with respect to the Killing field $X$ of the Hilbert form $\omega_1$.

\section{Functions $u(a)$ and $v(a)$ for spherically symmetric metrics}

The following notations and lemmas will be used in this section.
Let $F=|y|\phi\left(\frac{|x|^2}2,\, \frac{\langle x,\,y\rangle}{|y|}\right)$
be a spherically symmetric Finsler metric on $\mathbb{B}^2(\mu)$. Let
$$
r:=|y|,\quad t:=\frac{|x|^2}2,\quad s:=\frac{\langle x,\,y\rangle}{|y|},
\eqno(7.1)
$$
$$
r^i:=r_i:=\frac{y^i}{|y|},\quad x_i:=x^i, \quad s^i:=s_i:=x_i-sr_i.
\eqno(7.2)
$$
By a straightforward computation one obtains
$$
s_{y^i}=\frac{s_i}r
\eqno(7.3)
$$
where we have used (7.1) and (7.2).

\

{\bf Lemma 7.1}[9]\, {\em Let $f=f(r,\,t,\,s)$ be a function
on a domain $\mathcal{U}\subset \mathbb{R}^3$}. Then
$$
f_{x^i}=(r_i,\,s_i)\left(
\begin{array}{c}
  f_s+sf_t\\
  f_t
\end{array}
\right),\quad
f_{y^i}=(r_i,\,s_i)\left(
\begin{array}{c}
  f_r\\
  f_s/r
\end{array}
\right).
\eqno(7.4)
$$

\

{\bf Corollary 7.2}[9]\,\,{\em Let $F=|y|\phi\left(\frac{|x|^2}2,
\, \frac{\langle x,\,y\rangle}{|y|}\right)$ be a spherically symmetric
Finsler metric on $\mathbb{B}^2(\mu)$. Then
$$
F_{y^i}=\phi r_i+\phi_s s_i.
\eqno(7.5)
$$
and
$$
F_0=:F_{x^i}y^i=r^2\cdot(\phi_s+s\phi_t).
\eqno(7.6)
$$
}

The geodesic coefficients $G^i$ can be expressed by (cf [9],\,[15, Definition 3.3.8])
$$
G^i:=\frac{r^{2}}{2}(r^{i},s^{i})\left(\begin{array}{ccc}
\bar{u}\\ \bar{v}\\
\end{array}\right).
 \eqno(7.7)
$$
where
$$
\bar{v}:=\frac{s\phi_{ts}+\phi_{ss}-\phi_{t}}{\Delta},\quad \bar{u}
        =\frac{1}{\phi}[\phi_{s}+s\phi_{t}-(2t-s^{2})\phi_{s}\bar{v}].
\eqno(7.8)
$$
where
$$
\Delta=\phi-s\phi_{s}+(2t-s^{2})\phi_{ss}.\eqno(7.9)
$$
By a straightforward computation one obtains the following

\

{\bf Lemma 7.3}\,\,{\em Let $t$ and $s$ be functions satisfying (7.1). Then
$$
2t-s^2=\frac{(x^1y^2-x^2y^1)^2}{|y|^2}.
\eqno(7.10)
$$

}

Now we give an approach to calculate the normal forms
for known spherically symmetric Finsler metrics of constant flag curvature.

\

{\bf Step 1} \,\, First of all, let us calculate $a=a(t,\,s)$.
Most of known spherically symmetric  Finsler metrics of constant flag curvature
are projectively flat. Hence $\phi$ satisfies the following projectively flat equation [10]:
$s\phi_{ts}+\phi_{ss}-\phi_t=0.$
It follows that
$\phi-s\phi_s=h(2t-s^2).$
where $h$ is a function.  In this case, $a=a(z)$ where $z=2t-s^2$ (see (7.15) below).

\

Let $F=|y|\phi\left(\frac{|x|^2}2,\, \frac{\langle x,\,y\rangle}{|y|}\right)$
be a spherically symmetric Finsler metric on $\mathbb{B}^2(\mu)$.
We can express $x=(x^1,\,x^2)$ in the polar coordinate system,
$$
x^1=\rho\cos\theta,\qquad x^2=\rho\sin\theta.
\eqno(7.11)
$$
By a straightforward computation one obtains
$$
X:=\frac{\partial}{\partial\theta}=-\rho\sin\theta\frac{\partial}{\partial x^1}
+\rho\cos\theta\frac{\partial}{\partial x^2}=-x^2\frac{\partial}{\partial x^1}+x^1\frac{\partial}{\partial x^2}.
\eqno(7.12)
$$
From (7.11) we have
$$
F=\sqrt{p^2+q^2\rho^2}\phi\left(\frac{\rho}2,\,\frac{p\rho}{\sqrt{p^2+q^2\rho^2}}\right)
$$
It follows that
$$
X(F)=0.
\eqno(7.13)
$$
(7.13) tells us that $\frac{\partial}{\partial\theta}$
is a Killing field of $F$ [17] and all spherically symmetric Finsler surfaces
admit a non-zero Killing field $\frac{\partial}{\partial\theta}$.
Let $X$ be the vector field (7.2) on $\mathbb{B}^2(\mu)$.
By using (4.17) and (7.12), we have
$$
\hat{X}=-x^2\frac{\partial}{\partial x^1}+x^1\frac{\partial}{\partial x^2}-y^2\frac{\partial}{\partial y^1}+y^1\frac{\partial}{\partial y^2}.
\eqno(7.14)
$$
 From (2.1), (7.1), (7.2), (7.5), (7.10), (3.3) and (7.14), we obtain
$$
\begin{array}{ccl}
a_1
&=&
\left(F_{y^1}dx^1+F_{y^2}dx^2\right)\left(-x^2\frac{\partial}{\partial x^1}+x^1\frac{\partial}{\partial x^2}-y^2\frac{\partial}{\partial y^1}+y^1\frac{\partial}{\partial y^2}\right)\\
&=&
-F_{y^1}x^2+F_{y^2}x^1\\
&=&
-\left(\phi r_1+\phi_s s_1\right)x^2+\left(\phi r_2+\phi_s s_2\right)x^1=(\phi-s\phi_s)\frac{\sqrt{2t-s^2}}{|y|}.
\end{array}
\eqno(7.15)
$$

\

{\bf Step 2}\,\,
We now calculate $f(z)$ and $g(z)$ for fixed constant $K$ (cf (3.9) and (3.14)).
A direct calculation yields
[21, Proposition 3.2]
$$
D:=\det(g_{ij})=\phi^3\Delta
\eqno(7.16)
$$
where $\Delta$ is given in (3.9). Note that $\{\omega_1,\,\,\omega_2\}$
is the Berwald frame on the Finsler surface $\mathbb{B}^2(\mu)$.
It follows that $\omega_2=\frac{\sqrt D}F\left(-y^2dx^1+y^1dx^2\right)$ [1].
Together with (3.3), (7.14) and (7.16), we have
$$
a_2
=
\frac{\sqrt D}F(y^2x^2+y^1x^1)
=
\frac{\sqrt{\phi^3\Delta}}{\phi|y|}\langle x,\,y\rangle=s\phi^{\frac 12}\Delta^{\frac 12}.
\eqno(7.17)
$$
We express the geodesic coefficients $G^i$ by (cf [6])
$$
G^i=Py^i+Q^i
\eqno(7.18)
$$
where
$$
P=\frac{F_0}{2F}=\frac r{2\phi}\left(\phi_s+s\phi_t\right)=\frac r2(\bar{u}-s\bar{v}),
\eqno(7.19)
$$
$$
Q=\frac{r^2}2 \bar{v}x^i
\eqno(7.20)
$$
where we have used (7.2) and (7.7). By (7.18), we have the connection coefficients
$$
N_j^i:=\frac{\partial G^i}{\partial y^j}=P_{y^j}y^i+P\delta^i_j+x^i\left(\frac{r^2}2 \bar{v}\right)_{y^j}.
\eqno(7.21)
$$
On the other hand, we have (cf [2], Page 93)
$$
\omega_3=\frac{\sqrt D}{F^2}\left[y^1(dy^2+N_j^2dx^j)-y^2(dy^1+N_j^1dx^j)\right].
$$
Together with (3.3), (7.14) and (7.21), we get
$$
\begin{array}{ccl}
a_3
&=&
\frac{\sqrt D}{F^2}\left(|y|^2+x^2y^2N^1_1+x^1y^1N^2_2-x^2y^1N^2_1-x^1y^2N^1_2\right)
\\
&=&
\frac{\sqrt D}{F^2}\left(|y|^2+P\langle x,\,y\rangle +(I)\right).
\end{array}
\eqno(7.22)
$$
where
$$
\begin{array}{ccl}
(I):
&=&
x^1x^2y^2\left(\frac{r^2}2 \bar{v}\right)_{y^1}+x^1x^2y^1\left(\frac{r^2}2 \bar{v}\right)_{y^2}\\
&&
-(x^1)^2y^2\left(\frac{r^2}2 \bar{v}\right)_{y^2}-(x^2)^2y^1\left(\frac{r^2}2 \bar{v}\right)_{y^1}\\
&=&
-\bar{v}(x^1y^2-x^2y^1)^2+\frac{r^2}2\left(\bar{v}_{y^1}x^2-\bar{v}_{y^2}x^1\right)(x^1y^2-x^2y^1)\\
&=&
-\frac 12(2\bar{v}-s\bar{v}_s)(x^1y^2-x^2y^1)^2=-\frac{r^2}2(2\bar{v}-s\bar{v}_s)(2t-s^2)
\end{array}
\eqno(7.23)
$$
where we have used (7.10) and
$$
\bar{v}_{y^j}=\frac{\bar{v}_s}r\left(x^j-s\frac{y^j}r\right),
\qquad \bar{v}_{y^1}x^2-\bar{v}_{y^2}x^1=\frac{s\bar{v}_s}{r^2}(x^1y^2-x^2y^1).
\eqno(7.24)
$$
Plugging (7.19) and (7.23) into (7.22) yields
$$
\begin{array}{ccl}
a_3
&=&
\frac{\sqrt D}{F^2}\left[r^2+\frac{sr^2}2(\bar{u}-s\bar{v})-\frac{r^2}2(2\bar{v}-s\bar{v}_s)(2t-s^2)\right]
\\
&=&
\frac{\Delta^{\frac 12}}{2\phi^{\frac 12}}\left[2+s(\bar{u}-s\bar{v})-(2\bar{v}-s\bar{v}_s)(2t-s^2)\right]
\end{array}
\eqno(7.25)
$$
where  we have used (7.16). The main scalar $I$ of $F$ is given by
$$
\begin{array}{ccl}
I
&=&
\left[\frac FD\left(\frac{\partial}{\partial y^1}\log\sqrt D\right)dx^1+\frac FD\left(\frac{\partial}{\partial y^2}\log\sqrt D\right)dx^2\right]
\left(-F_{y^2}\frac{\partial}{\partial x^1}+F_{y^1}\frac{\partial}{\partial x^2}\right)\\
&=&
\frac F{2D^{\frac 32}}(F_{y^1}D_{y^2}-F_{y^2}D_{y^1})\\
&=&
\frac{\phi^2D_s}{2D^{\frac 32}}(r_2s_1-r_1s_2)=\sqrt{2t-s^2}\frac{\phi^2D_s}{2D^{\frac 32}}
\end{array}
\eqno(7.26)
$$
where we have made use of (7.5) and the following fact:
$$
D_{y^i}=D_s\frac{s_i}r.
$$
Together with (7.25), we have
$$
a_3I=\frac{\Delta^{\frac 12}}{2\phi^{\frac 12}}T\sqrt{2t-s^2}\frac{\phi^2D_s}{2D^{\frac 32}}
    =\frac{\Psi}{4\Delta\phi}\sqrt{2t-s^2}T
\eqno(7.27)
$$
where
$$
T:=2+s(\bar{u}-s\bar{v})-(2\bar{v}-s\bar{v}_s)(2t-s^2)
$$
and
$$
\Psi:=3\phi_s\Delta+\phi\Delta_s=\frac{D_s}{\phi^2}.
$$
By Lemma 7.1, we have
$$
I_{x^j}y^j=r(I_s+sI_t).
\eqno(7.28)
$$
 Note that $I$ is homogeneous of degree zero with respect to $y$.
Hence $y^j\frac{\partial I}{\partial y^j}=0.$
Together with (2.7), (7.21) and (7.28), we get
$$
J=e_1(I)=\frac 1{\phi}(I_s+sI_t)-\frac{(II)}F
\eqno(7.29)
$$
where
$$
e_1:=\frac 1Fy^j\left(\frac{\partial}{\partial x^j}-N^k_j\frac{\partial}{\partial y^k}\right)
$$
and
$$
\begin{array}{ccl}
(II):
&=&
\frac{x^ky^j}2(r^2\bar{v})_{y^j}\frac{\partial I}{\partial y^k}\\
&=&
r^2x^k\left(\bar{v}+\frac{y^j}2\bar{v}_{y^j}\right)\frac{\partial I}{\partial y^k}
=r\bar{v}I_s(2t-s^2)
\end{array}
\eqno(7.30)
$$
where we have used
$$
y^j\bar{v}_{y^j}=0,\qquad \frac{\partial I}{\partial y^j}=I_s\frac{s_j}r.
$$
Substituting (7.30) into (7.29) yields
$$
J=\frac 1{\phi}\square I
\eqno(7.31)
$$
where
$$
\square:=s\frac{\partial}{\partial t}+\left[1-(2t-s^2)\bar{v}\right]\frac{\partial}{\partial s}.
$$
By a straightforward computation one obtains
$$
\square\left(\sqrt{2t-s^2}\right)=s\sqrt{2t-s^2}\bar{v}.
\eqno(7.32)
$$
Substituting (7.26) into (7.31) and using (7.17) and (7.32) we have
$$
\begin{array}{ccl}
a_2J
&=&
s\phi^{\frac 12}\Delta^{\frac 12}\cdot \frac 1{2\phi}
\square\left(\frac{\sqrt{2t-s^2}\Psi}{\phi^{\frac 12}\Delta^{\frac 32}}\right)\\
&=&
\frac{s\sqrt{2t-s^2}}{4\phi\Delta^2}
\left[2\Delta(\square\Psi+s\Psi\bar{v})-\Psi(\Delta\square\log\phi+3\square\Delta)\right].
\end{array}
\eqno(7.33)
$$
By using (7.17),(7.25), (7.27) and (7.33), we obtain
$$
f=f(z),\qquad g=g(z).
$$

{\bf Step 3}\,\, Let $u(a)^2:=f(z(a))$ and $v(a)=\frac{g(z(a))}{f(z(a))}$.
Substituting these these into (4.16), (5.10) and (6.15), we obtain
the normal forms for known spherically symmetric Finsler metrics of constant flag curvature.

For example, for the Funk metric on the unit disk $\mathbb{D}^2$,
we find $u(a)=\sqrt{1+4a^2}$ and $v(a)=\frac{-3a}{1+4a^2}$.


\end{document}